\documentclass[12pt]{amsart}

\usepackage{amssymb}
\usepackage{amscd}

\newcommand{\bt}{\begin{theorem}}
\newcommand{\et}{\end{theorem}}
\newcommand{\bp}{\begin{proposition}}
\newcommand{\ep}{\end{proposition}}
\newcommand{\bq}{\begin{question}}
\newcommand{\eq}{\end{question}}
\newcommand{\bl}{\begin{lemma}}
\newcommand{\el}{\end{lemma}}
\newcommand{\br}{\begin{result}}
\newcommand{\er}{\end{result}}
\newcommand{\be}{\begin{equation}}
\newcommand{\ee}{\end{equation}}
\newcommand{\bc}{\begin{corollary}}
\newcommand{\ec}{\end{corollary}}
\newcommand{\bex}{\begin{example}}
\newcommand{\eex}{\end{example}}

\newtheorem{theorem}{Theorem}[section]
\newtheorem{corollary}[theorem]{Corollary}
\newtheorem{lemma}[theorem]{Lemma}
\newtheorem{proposition}[theorem]{Proposition}
\newtheorem{result}[theorem]{Result}
\newtheorem{example}[theorem]{Example}
\newtheorem{question}[theorem]{Question}

\numberwithin{equation}{section}

\newcommand{\N}{\mathbb{N}}

\newcommand{\Z}{\mathbb{Z}}

\newcommand{\cB}{\mathcal{B}}

\newcommand{\cH}{\mathcal{H}}

\newcommand{\cL}{\mathcal{L}}

\newcommand{\cR}{\mathcal{R}}

\newcommand{\e}{\varepsilon}
\newcommand{\f}{\varphi}
\newcommand{\om}{\omega}

\newcommand{\epr}{\hspace{\fill}$\Box$}

\newcommand{\med}{\medskip}
\newcommand{\sm}{\smallskip}

\begin{document}
\noindent {\em Semigroup Forum} (2015) {\bf 91}, 117 -- 127\\
DOI 10.1007/s00233-014-9649-1
\vspace{0.02in}\\
{\em arXiv version}: layout, fonts, pagination and numbering of sections, lemmas, theorems and formulas may vary from the SF published version
\vspace{0.04in}\\
\title[Perfect congruences on bisimple $\om$-semigroups]{Perfect congruences on bisimple $\om$-semigroups*}\thanks{*This research was partially supported by NSF grant DMS-1156612}
\author[S. M. Goberstein, K. Grimshaw, A. Kling, T. Landry, F. Li ]{Simon M. Goberstein, Katherine Grimshaw, Anthony Kling, \vspace{0.05in}\\Therese Landry, Freda Li}
\address{Department of Mathematics and Statistics,
California State University, Chico, CA 95929, e-mail: sgoberstein@csuchico.edu}
\begin{abstract}
A congruence $\e$ on a semigroup $S$ is perfect if for any congruence classes  $x\e$ and $y\e$ their product as subsets of $S$ coincides (as a set) with the congruence class $(xy)\e$. Perfect congruences on the bicyclic semigroup were found in \cite{key7}. Using the structure of bisimple $\om$-semigroups determined in \cite{key25} and the description of congruences on these semigroups found in \cite{key20} and \cite{key1}, we obtain a complete characterization of perfect congruences on all bisimple $\om$-semigroups, substantially generalizing the above mentioned result of \cite{key7}.  
\vspace{0.1in}
\vspace{0.05in}
\end{abstract}
\maketitle
\font\caps=cmcsc10 scaled \magstep1   
\def\bfseries{\normalsize\caps}
\vspace{-0.2in}
\section{Introduction}
\sm 
Let $S$ be a semigroup. Recall that for all $A,B\subseteq S$, the {\em set product} of $A$ and $B$ is defined by the formula $AB=\{ab:a\in A, b\in B\}$ (see \cite[\S 1.1]{key3}).  Let $\e$ be a congruence on $S$. The $\e$-class containing an element $x$ of $S$ will be denoted by $x\e$, so $x\e=\{x'\in S: (x,x')\in\e\}$.  Since the inclusion $(x\e)(y\e)\subseteq(xy)\e$ holds for all $x,y\in S$, it is immediate (and well known) that one can unambiguously define an operation $\ast$ on the set $S/\e$ of all $\e$-classes by the formula $(x\e)\ast(y\e)=(xy)\e$, and with respect to this operation  $S/\e$ becomes a semigroup called the {\em factor semigroup} of $S$ modulo $\e$.  It must be emphasized that, in general, the set product $(x\e)(y\e)$ may be properly contained in $(xy)\e$ for some $x,y\in S$. A congruence $\e$ on $S$ is {\em perfect} if $(x\e)(y\e)=(xy)\e$, that is, if $(x\e)(y\e)=(x\e)\ast(y\e)$ for all $x,y\in S$, and $S$ is said to be a {\bf\em perfect semigroup} if all congruences on $S$ are perfect. In this paper, we will make a clear distinction between $(x\e)\ast(y\e)$, the result of the operation of the factor semigroup $S/\e$ applied to $\e$-classes $x\e$ and $y\e$, and the set product $(x\e)(y\e)$ of subsets $x\e$ and $y\e$ of $S$.  

The concept of a `perfect congruence' was originally introduced by Wagner \cite{key30} for partial groupoids. Perfect congruences on semigroups and on some other algebraic structures have been studied by a number of authors. Clearly, groups are perfect semigroups, so `perfectness' can be viewed as a `group-like' property of a semigroup. Since inverse semigroups represent one of the most important generalizations of groups, it is natural to consider the problem of characterizing perfect congruences on inverse semigroups and identifying those inverse semigroups which are perfect. For instance, the structure of perfect Clifford semigroups was described by Fortunatov \cite{key5} and that of perfect finite inverse semigroups by Hamilton and Tamura \cite{key12}. Both of these results were generalized by Goberstein \cite{key10} who determined the structure of perfect completely semisimple inverse semigroups. 

The  {\em bicyclic semigroup} $\cB(a,b)$ is a semigroup with identity $1$ generated by the two-element set $\{a,b\}$ and given by one defining relation $ab=1$ (see \cite[\S 1.12]{key3}). It is well known that $\cB(a,b)$ is a bisimple inverse semigroup with no nontrivial subgroups, each of its elements has a unique representation in the form $b^ma^n$ where $m$ and $n$ are nonnegative integers (and $a^0=b^0=1$), the semilattice of idempotents of $\cB(a,b)$ is a chain: $1>ba>b^2a^2>\cdots$, and if $\e$ is an arbitrary congruence on $\cB(a,b)$, then either $\e$ is the equality relation on $\cB(a,b)$ or $\cB(a,b)/\e$ is a cyclic group (see \cite[Lemma 1.31, Corollary 1.32, and Theorem 2.53 and its proof]{key3}). It follows from \cite[Theorem 2.54]{key3} that an inverse semigroup is completely semisimple if and only if it has no subsemigroup isomorphic to the bicyclic semigroup. Thus, in view of \cite{key10}, it is appropriate to ask: Which congruences on the bicyclic semigroup are perfect? The answer to this question is known -- as shown in \cite{key7}, the only congruence on the bicyclic semigroup which is not perfect is the minimum group congruence. The next natural step is to consider the problem of characterizing perfect congruences on arbitrary bisimple inverse semigroups whose idempotents form a descending chain isomorphic to the chain of idempotents of the bicyclic semigroup.  Such semigroups are called {\em  bisimple $\om$-semigroups}, and the purpose of this article is to identify and describe perfect congruences on bisimple $\om$-semigroups.

In Section 2 we discuss the connection between divisibility relations and perfect congruences on semigroups and gather basic information about bisimple $\om$-semigroups and their congruences. Our new results are contained in Sections 3 and 4. Let $S$ be a bisimple $\om$-semigroup. We prove that each idempotent-separating congruence on $S$ is perfect (Theorem \ref{302}) and determine which group congruences on $S$ are perfect (Proposition \ref{408} and Theorem \ref{410}). Since every congruence on $S$ is either a group congruence or an idempotent-separating one  \cite[Theorem 1.3]{key20}, our results provide a complete characterization of perfect congruences on bisimple $\om$-semigroups.  

We use \cite{key3} and \cite{key15} as standard references for the algebraic theory of semigroups and refer to \cite{key23} for an extensive treatment of the theory of inverse semigroups.
\sm
\section{Preliminaries}
\sm 
We begin by recalling a few common notational conventions which will be used throughout the paper. Let $\rho$ be an arbitrary binary relation on some set $U$. If $X$ is a subset of $U$, then $X\rho=\{u\in U:(\exists\,x\in X)\;(x,u)\in\rho\}$, and if $X=\{x\}$, we will write $x\rho$ instead of $\{x\}\rho$. As usual, $\rho^{-1}=\{(v,u)\in U\times U: (u,v)\in\rho\}$. Therefore for any subset $Y$ of $U$, we have $Y\rho^{-1}=\{u\in U:(\exists\,y\in Y)\;(u,y)\in\rho\}$; in particular, if $\f:U\to U$ is an arbitrary mapping, then $Y\f^{-1}=\{u\in U:u\f\in Y\}$. The identity mapping on $U$ will be denoted by $1_U$ (we will not distinguish $1_U$ from the equality relation $\{(u, u) : u\in U\}$ on $U$).
\vspace{0.03in}\\
\indent Let $S$ be a semigroup. Recall that if $x=yz$ for some $x, y, z\in S$, it is common to say that $y$ is a {\em left} and $z$ a {\em right divisor} of $x$ while $x$ is a {\em right multiple} of $y$ and a {\em left multiple} of $z$. Let 
\[\delta_l=\{(a,b)\in S\times S:b\in Sa\}\text{ and }\delta_r=\{(a,b)\in S\times S : b\in aS\},\] 
so $(a,b)\in\delta_l$ if and only if $b$ is a left multiple of $a$, and $(a,b)\in\delta_r$ if and only if $b$ is a right multiple of $a$. We refer to $\delta_l$ and $\delta_r$ as the {\em left} and {\em right divisibility relations} on $S$. In the notation of the preceding paragraph, $a\delta_l=Sa$ and $a\delta_r=aS$ for all $a\in S$. It is clear that if $\f$ is a homomorphism of $S$ to a semigroup $T$ and if $(a,b)\in\delta_l$ [$(a,b)\in\delta_r$] in $S$, then $(a\f,b\f)\in\delta_l$ [$(a\f,b\f)\in\delta_r$] in $T$.  Note also that if $S$ is regular (or if $S$ contains the identity element), then for all $a,b\in S$, we have $(a,b)\in\delta_l$  if and only if $Sb\subseteq Sa$, and $(a,b)\in\delta_r$  if and only if $bS\subseteq aS$, so that $\delta_l\cap\delta^{-1}_l=\cL$ and $\delta_r\cap\delta^{-1}_r=\cR$.

As observed in \cite{key7}, there is a connection between divisibility relations on semigroups and perfect congruences of certain types. Let $S$ be a semigroup. Recall that $\e$ is a {\em left} [{\em right}] {\em cancellative} congruence or a {\em group} congruence on $S$ if $S/\e$ is a left [right] cancellative semigroup or a group, respectively (see \cite[Proposition 1.7]{key3} for a general description of such terminology). 

From the proof of \cite[Lemma 9]{key7}, it can be deduced that if $\e$ is an arbitrary perfect congruence on a semigroup $S$, then $\delta_r\circ\e=\e\circ\delta_r$ and $\delta_l\circ\e=\e\circ\delta_l$.  Using this observation, we can state Lemma 9 and Corollary 10 of \cite{key7} as follows:
\br\label{201}{\rm\cite[Lemma 9]{key7}} Let $\e$ be a left {\rm[}right{\rm]} cancellative congruence on a semigroup $S$. Then $\e$ is perfect if and only if $\delta_r\circ\e=\e\circ\delta_r$ {\rm[}$\delta_l\circ\e=\e\circ\delta_l${\rm]}. 
\er
\br\label{202}{\rm\cite[Corollary 10]{key7}} Let $\e$ be a group congruence on a semigroup $S$. Then $\e$ is perfect if and only if $\e\circ\delta_r=S\times S$ {\rm[}$\delta_l\circ\e=S\times S${\rm]}. 
\er
{\bf Remark 1.} The right divisibility relation on a semigroup $S$ is denoted in \cite{key7}  by $\tau_d$ (instead of $\delta_r$), and Lemma 9 of \cite{key7}  states that a {\em right} cancellative congruence $\e$ on $S$ is perfect if and only if $\tau_d\circ\e=\e\circ\tau_d$ while Corollary 10 asserts that a group congruence $\e$ on $S$ is perfect if and only if $\tau_d\circ\e=S\times S$. The discrepancy between these statements and those of Results \ref{201} and \ref{202}, respectively, is due to the fact that in the composition $\alpha\circ\beta$ of binary relations $\alpha$ and $\beta$  on $S$ we view $\alpha$ as the first factor whereas in \cite{key7} the first factor is considered to be $\beta$.
\vspace{0.02in}\\ 
\indent Denote by $\N$ the set of all nonnegative integers. For all $(m,n),(p,q)\in \N\times\N$, define 
\[(m,n)(p,q)=(m+p-r,n+q-r),\]
where $r=\min\{n,p\}$. Then $\N\times\N$ endowed with this multiplication is a semigroup which we denote by $B$. Since $B$ is isomorphic to the bicyclic semigroup $\cB(a,b)$ (see \cite[page 45]{key3}), we will identify $B$ with the bicyclic semigroup. It is easily seen  that for $(m,n), (i,j)\in B$, we have $((m,n), (i,j))\in\delta_r$ [$((m,n), (i,j))\in\delta_l$] if and only if $m\leq i$ [$n\leq j$]; this simple observation will be used below without mention.
\vspace{0.02in}\\
\indent Let $G$ be a group and $\alpha$ an endomorphism of $G$. For all $(m,g,n), (p,h,q)\in\N\times G\times\N$, set 
\setcounter{equation}{2} 
\be(m,g,n)(p,h,q)=(m+p-r,(g\alpha^{p-r})( h\alpha^{n-r}),n+q-r),
\label{eq:1st}
\ee
where $r=\min\{n,p\}$ and $\alpha^0$ stands for the identity automorphism $1_G$ of $G$. As in \cite{key15}, we denote the set $\N\times G\times\N$ equipped with multiplication defined by (\ref{eq:1st}) by $BR(G,\alpha)$ and call it the {\em Bruck-Reilly extension of $G$ determined by $\alpha$}. One can readily check that $\{(n, e_G, n):n\in\N\}$ is the set of idempotents of $BR(G,\alpha)$ and $(0, e_G, 0)>(1, e_G, 1)>\cdots>(n, e_G, n)>\cdots$, so the idempotents of $BR(G,\alpha)$ form a semilattice isomorphic to the semilattice of idempotents of the bicyclic semigroup (here and elsewhere in the paper $e_G$ denotes the identity element of the group $G$). It is easily seen that $BR(G,\alpha)$ is an inverse semigroup with identity $(0, e_G, 0)$ and $(m,g,n)^{-1}=(n,g^{-1},m)$ for all $(m,g,n)\in BR(G,\alpha)$. Define a mapping $\psi: BR(G,\alpha)\to B$ by the rule: $(m,g,n)\psi=(m,n)$ for all $(m,g,n)\in BR(G,\alpha)$. Then $\psi$ is a homomorphism of $BR(G,\alpha)$ onto $B$ which we will call the {\em group-forgetful} homomorphism of $BR(G,\alpha)$ onto $B$ (of course, if $G=\{e_G\}$, then $\psi$ is injective and so is an isomorphism of $BR(G,\alpha)$ onto $B$).   
\setcounter{theorem}{3}
\br\label{204}{\rm\cite[Theorems 2.2 and 3.5]{key25}} If $G$ is a group and $\alpha$ an endomorphism of $G$, then $BR(G,\alpha)$ is a bisimple $\om$-semigroup. Conversely, if $S$ is any bisimple $\om$-semigroup and $G$ is its group of units, there is an endomorphism $\alpha$ of $G$ such that $S$ is isomorphic to $BR(G,\alpha)$.
\er
Let $S$ be an inverse semigroup. Denote by $E_S$ the semilattice of idempotents of $S$. Clearly, if $\e$ is a group congruence on $S$, all idempotents of $S$ are contained in a single $\e$-class. Let $\sigma_S=\{(x,y)\in S\times S: (\exists\,e\in E_S)\;ex=ey\}$. By \cite[Theorem V.3.1]{key15}, $\sigma_S$ is the minimum group congruence on $S$ (that is, $\sigma_S\subseteq\e$ for any group congruence $\e$ on $S$). Recall that a congruence $\e$ on $S$ is  {\em idempotent-separating} if $(e,f)\in\e$ implies $e=f$ for all $e,f\in E_S$. Of course, if $S$ is not a group, no group congruence on $S$ can be idempotent-separating. By \cite[Theorem V.3.2]{key15}, there exists a maximum idempotent-separating congruence on $S$ -- it is the largest congruence on $S$ contained in the Green's relation $\cH$. According to \cite[Lemma 1.2]{key20}, if $S=BR(G,\alpha)$, then $\cH$ is a congruence on $S$ and $S/\cH\cong B$. 
\vspace{0.02in}\\
\indent Let $G$ be a group, $\alpha$ an endomorphism of $G$, and $H$ a subgroup of $G$. Recall that $H$ is said to be $\alpha$-{\em admissible} if $H\alpha\subseteq H$ (see, for instance, \cite[Section 2.4]{key13}). As in \cite{key1}, we will also say that $H$ is {\em$\alpha$-invariant} if $H\alpha^{-1}=H$; clearly, an $\alpha$-invariant subgroup $H$ of $G$ is $\alpha$-admissible since $H\alpha^{-1}=H$ implies $H\alpha=(H\alpha^{-1})\alpha\subseteq H$.  For each congruence $\e$ on $BR(G,\alpha)$, define
\setcounter{equation}{4} 
\be N_{\e}=\{a\in G: ((0, a, 0), (0, e_G, 0))\in\e\},
\label{eq:2nd}
\ee
and observe that, according to \cite[Lemma 2.1]{key20}, $N_{\e}$  is an $\alpha$-admissible normal subgroup of $G$.
\setcounter{theorem}{5}
\br\label{205}{\rm\cite[Theorem 1.3]{key20}} Let $G$ be a group and $\alpha$ an endomorphism of $G$. Then every congruence on $BR(G,\alpha)$ is either a group congruence or an idempotent-separating one.
\er
In view of Result \ref{205}, the problem of characterizing perfect congruences on bisimple $\om$-semigroups can be solved by considering it separately for idempotent-separating congruences and for group congruences on such semigroups. This will be done in the next two sections. 
\sm
\section{Idempotent-separating congruences on bisimple $\om$-semigroups are perfect}
\sm
Let $G$ be a group and $\alpha$ an endomorphism of $G$. There is a close connection between $\alpha$-admissible normal subgroups of $G$ and idempotent-separating congruences on $BR(G,\alpha)$.

\br\label{301}{\rm\cite[Lemma 2.3]{key20}} Let $G$ be a group, $\alpha\in{\rm End}(G)$, and $S=BR(G,\alpha)$. 
\vspace{0.02in}\\
{\rm(i)} Let $\e$ be an idempotent-separating congruence on $S$. Then 
\[ ((m, g, n), (p, h, q))\in\e\Longleftrightarrow m=p,\; n=q\,\text{  and  }\,gh^{-1}\in N_{\e}.\]
{\rm(ii)} For any $\alpha$-admissible normal subgroup $N$ of $G$ there exists an idempotent-separating congruence $\e$ on $S$ such that $N=N_{\e}$.
\er
According to statement (i) of Result \ref{301}, if $\e$ is an idempotent-separating congruence on $BR(G,\alpha)$ and if $x=(m,g,n)\in BR(G,\alpha)$, then $x\e=\{(m,ag,n):a\in N_{\e}\}$. Moreover, using statement (i) of Result \ref{301}, one can make statement (ii) a bit more precise. Namely, for any $\alpha$-admissible normal subgroup $N$ of $G$ there is a unique idempotent-separating congruence $\e_{\!_N}$ on $S$ such that $N=N_{\e_{\!_N}}$; this congruence $\e_{\!_N}$  is defined as follows: 
\[ ((m, g, n), (p, h, q))\in\e_{\!_N}\Longleftrightarrow m=p,\; n=q\,\text{  and  }\,gh^{-1}\in N.\]
Thus if $\e$ is an idempotent-separating congruence on $BR(G,\alpha)$, then $\e_{\!_{N_{\e}}}=\e$, and if $N$ is an $\alpha$-admissible normal subgroup of $G$, then $N_{\e_{\!_N}}=N$.

\bt\label{302} Let $G$ be a group and $\alpha$ an endomorphism of $G$. Every idempotent-separating congruence on $BR(G,\alpha)$ is perfect. 
\et
{\bf Proof}. Let $S=BR(G,\alpha)$, and  let $\e$ be any idempotent-separating congruence on $S$. Denote, for short, $N=N_{\e}$. Take arbitrary $x,y\in S$. Thus $x=(m,g,n)$ and $y=(p,h,q)$ for some $m,n,p,q\in\N$ and $g,h\in G$. As noted in the paragraph following Result \ref{301},  $x\e=\{(m,ag,n):a\in N\}$ and $y\e=\{(p,bh,q):b\in N\}$. By definition of multiplication in $S$,
\setcounter{equation}{2} 
\be xy=(m+p-r, (g\alpha^{p-r})(h\alpha^{n-r}),n+q-r),
\label{eq:3rd}
\ee
\be (x\e)(y\e)=\{(m+p-r, (ag)\alpha^{p-r}\cdot(bh)\alpha^{n-r},n+q-r): a, b\in N\}
\label{eq:4th}
\ee
and also
\be (xy)\e=\{(m+p-r, c(g\alpha^{p-r})(h\alpha^{n-r}),n+q-r):c\in N\},
\label{eq:5th}
\ee
 where $r={\rm min}\{n,p\}$. Our goal is to prove that $(xy)\e\subseteq(x\e)(y\e)$. From (\ref{eq:4th}) and (\ref{eq:5th}) it follows that $(xy)\e\subseteq(x\e)(y\e)$ if and only if for all $c\in N$ there exist $a,b\in N$ such that 
\be
  c(g\alpha^{p-r})(h\alpha^{n-r})=(ag)\alpha^{p-r}\cdot(bh)\alpha^{n-r},
\label{eq:6th}  
\ee
and since $\alpha$ is an endomorphism of $G$, (\ref{eq:6th}) is equivalent to 
\be 
c(g\alpha^{p-r})=(a\alpha^{p-r})(g\alpha^{p-r})(b\alpha^{n-r}).
\label{eq:7th}
\ee
Take an arbitrary $c\in N$. Let us show that there exist $a,b\in N$ such that (\ref{eq:7th}) holds. 
\vspace{0.03in}\\
\indent {\bf Case 1:} $p\leq n$.
\vspace{0.02in}\\
\indent In this case, $r=p$ and so $\alpha^{p-r}=\alpha^0$. Recall that $\alpha^0$ is the identity automorphism $1_G$ of $G$. Therefore (\ref{eq:7th}) becomes $cg=ag(b\alpha^{n-r})$, which holds if we take $a=c$ and $b=e_G$.
\vspace{0.03in}\\
\indent {\bf Case 2:} $n<p$.
\vspace{0.02in}\\
\indent Here $r=n$ and hence $\alpha^{n-r}=\alpha^0=1_G$. Therefore (\ref{eq:7th}) becomes
\be
c(g\alpha^{p-r})=(a\alpha^{p-r})(g\alpha^{p-r})b,
\label{eq:8th}
\ee
where $p-r>0$. Note that $(g\alpha^{p-r})^{-1}c(g\alpha^{p-r})\in N$ since $N$ is a normal subgroup of $G$. It follows that (\ref{eq:8th}) holds if we choose $b=(g\alpha^{p-r})^{-1}c(g\alpha^{p-r})$ and $a=e_G$.
\vspace{0.02in}\\
\indent We have shown that $(xy)\e\subseteq(x\e)(y\e)$ for all $x,y\in S$. Thus $\e$ is a perfect congruence on $S$. The proof is complete. \epr
\sm
\section{Perfect group congruences on bisimple $\om$-semigroups}
\sm
Since a homomorphic image of $B$ is either isomorphic to $B$ or is a cyclic group (see \cite[Corollary 1.32]{key3}), all congruences on $B$ different from $1_B$ are group congruences, and they can be readily described as follows. Let $k\in\N$. For $(m,n), (p,q)\in B$, define $((m,n), (p,q))\in\zeta_0$ if and only if $n-m=q-p$, and for $k\geq 1$, set $((m,n), (p,q))\in\zeta_k$ if and only if $n-m\equiv q-p\,({\rm mod}\,k)$. Then $\zeta_0=\sigma_B$ and $B/\zeta_0$ is an infinite cyclic group, and if $k\geq 1$, then $B/\zeta_k$ is a finite cyclic group of order $k$. Thus $\{\zeta_k: k\in\N\}$ is the set of all congruences on $B$ different from $1_B$. By \cite[Theorem 11]{key7}, $\zeta_k$ is a perfect congruence for each $k\geq 1$, but $\zeta_0$ is not perfect. 

Let $S$ be an inverse semigroup. An inverse subsemigroup $U$ of $S$ is called {\em normal} if $E_S\subseteq U$ and  if $x^{-1}Ux\subseteq U$ for all $x\in S$  (cf. \cite[Definition III.1.3]{key23}). If for all $u,x\in S$, from $u,ux\in U$ it follows that $x\in U$, it is common to say that $U$ is a {\em unitary} inverse subsemigroup of $S$. 

Let $S$ be an inverse semigroup and $U$ a normal unitary inverse subsemigroup of $S$. As noted in \cite{key1}, it follows from \cite[\S 7.4]{key4} (see also \cite[Theorem III.1.8]{key23}) that by setting $(x,y)\in\e_{U}$ if and only if $xy^{-1}\in U$ for $x,y\in S$, we obtain a group congruence $\e_{U}$ on $S$ such that $U$ is the identity element of $S/\e_{U}$, and conversely, if $\e$ is a group congruence on $S$ and $U$ is the identity element of $S/\e$, then $U$ is a normal unitary inverse subsemigroup of $S$, and hence $\e=\e_{U}$.

In what follows we will use a description of group congruences on bisimple $\om$-semigroups obtained in \cite{key1}; for convenience of reference, we reproduce it here (in a slightly modified form).
\br\label{401}{\rm(from \cite[Theorem 1 and its proof]{key1})} Let $G$ be a group, $\alpha$ an endomorphism of $G$, and $S=BR(G,\alpha)$. Suppose that an $\alpha$-invariant normal subgroup $N$ of $G$, an element $z$ of $G$, and a nonnegative integer $k$ are such that the following conditions hold: 
\setcounter{equation}{1}
\be 
N(z\alpha)=Nz,
\label{eq:9th}
\ee
and 
\be 
(\forall\,g\in G)\,g^{-1}(Nz)(g\alpha^k)=Nz.
\label{eq:10th} 
\ee
Define $U(N,z,k)\subseteq S$ as follows: if $k\geq 1$, then
\be U(N,z,k)=\left\{(m,g,n)\in S:m\equiv n\,({\rm mod}\,k)\text{ and } g\in Nz^l \text{ where }l=(n-m)/k\right\},
\label{eq:11th}
\ee
and if $k=0$, then
\be U(N,z,0)=\{(m,g,m)\in S: g\in N\}.
\label{eq:12th}
\ee
Then $U(N,z,k)$ is a normal unitary inverse subsemigroup of $S$, and therefore $\e_{U(N,z,k)}$ is a group congruence on $S$.
\vspace{0.01in}\\
\indent Conversely, let $\e$ be a group congruence on $S$. Then the $\alpha$-admissible normal subgroup $N_{\e}$ of $G$ given by formula {\rm(\ref{eq:2nd})} is $\alpha$-invariant. Let $U_{\e}$ be the identity element of $S/\e$.  Consider $M=\{m\in\N:(0,g,m)\in U_{\e}\text{ for some }g\in G\}$, and define $k_{\e}\in\N$ as follows: if $M=\{0\}$, then $k_{\e}=0$, and if $M\ne\{0\}$, set $k_{\e}=\min (M\setminus\{0\})$. Finally, fix some $z_{\e}\in G$ such that $(0,z_{\e},k_{\e})\in U_{\e}$. Then for $N=N_{\e}$, $z=z_{\e}$, and $k=k_{\e}$, formulas {\rm(\ref{eq:9th})} and {\rm(\ref{eq:10th})} hold, and  $U_{\e}$ coincides with $U(N,z,k)$ defined by {\rm(\ref{eq:11th})} and {\rm(\ref{eq:12th})}, so that $\e=\e_{U(N,z,k)}$. 
\er
Let $G$ be a group, $\alpha\in{\rm End}(G)$, and $S=BR(G,\alpha)$, and let $N$, $z$, $k$ be as in Result \ref{401}. To slightly shorten the notation, we will write $\e_{(N,z,k)}$ instead of  $\e_{U(N,z,k)}$. By Result \ref{401}, if $\e$ is a group congruence on $S$, there exist unique $N$ and $k$ such that $\e=\e_{(N,z,k)}$ for some $z\in G$; although $z$ here, in general, is not unique, it is easy to check that if we also have $\e=\e_{(N,z',k)}$ for some $z'\in G$, then $z'\in Nz$, that is, the coset $Nz$ is determined by $\e$ uniquely (this fact is also an immediate corollary of \cite[Theorem 2]{key1}). Thus the formula $\e=\e_{(N,z,k)}$ establishes a one-to-one correspondence between the set of group congruences $\e$ on $S$ and the set of ordered triples $(N, Nz, k)$ where $N$ is an $\alpha$-invariant subgroup of $G$, $z\in G$, and $k\in\N$ such that {\rm(\ref{eq:9th})} and {\rm(\ref{eq:10th})} hold. It is easily seen (and also follows from \cite[Lemma on page 355]{key1}) that $N(z^l)\alpha^n=Nz^l$ for all $l\geq 1$ and $n\geq 0$.  {\em The notation and observations of this paragraph will be used without any additional explanation through the rest of this section.}
\vspace{0.02in}\\
\indent 
It is not difficult to deduce from Result \ref{401} (see \cite[Corollary]{key1}) that for all $(m,g,n)\in S$, the $\e_{(N,z,k)}$-class of $(m,g,n)$ can be described as follows: if $k\geq 1$, then 
\be
(m,g,n)\e_{(N,z,k)}\!=\!\left\{(p,h,q)\in S: (\exists\,l\!\in\!\Z)\left[q-p=n-m+kl\text{ and }h\alpha^n\!\in\!Nz^l(g\alpha^q)\right]\right\},
\label{eq:13th}
\ee
and if $k=0$, then
\be (m,g,n)\e_{(N,z,0)}\!=\!\{(p,h,q)\in S: q-p=n-m\text{ and }h\alpha^n\!\in\! N(g\alpha^q)\}.
\label{eq:14th}
\ee
\setcounter{theorem}{7}
\bp\label{408} No congruence $\e_{(N,z,0)}$ on $S$ is perfect. In particular, the minimum group congruence $\sigma_S$ on $S$ is not perfect.
\ep
{\bf Proof.} Let $\e=\e_{(N,z,0)}$. Suppose that $\e$ is perfect. Then, by Result \ref{202}, $\e\circ\delta_r=S\times S$. It follows, in particular, that $((1,e_G, 0), (0, e_G, 0))\in\e\circ\delta_r$. Thus there is $(p, h, q)\in S$ such that $((1,e_G, 0), (p, h, q))\in\e$ and $((p, h, q), (0, e_G, 0))\in\delta_r$. In view of (\ref{eq:14th}), $p=q+1\geq 1$. Using the group-forgetful homomorphism of $S$ onto $B$, from $((p, h, q), (0, e_G, 0))\in\delta_r$ we deduce that $((p, q), (0, 0))\in\delta_r$ in $B$. Then $p\leq 0$ whence $p=0\ngeq 1$. This contradiction shows that $\e$ is not perfect. By \cite[Theorem 4]{key1}, $\sigma_S=\e_{(N_0,e_G,0)}$ where $N_0=\{g\in G:g\alpha^n=e_G\text{ for some }n\geq 1\}$. Thus $\sigma_S$ is not perfect.
\epr
\vspace{0.02in}\\
\indent In view of Proposition \ref{408}, perfect group congruences on $S$ can be found only among congruences $\e_{(N,z,k)}$ such that $k\geq 1$. 
\bl\label{409}
Let $\e=\e_{(N,z,k)}$ with $k\geq 1$. Suppose that $\e$ is perfect. Then for any $(i,f,j)\in S$ there exist $m,n\in\N$ and $g,h\in G$ such that $(0,h,n)\in(i,f,j)\e$ and $(m,g,0)\in(i,f,j)\e$.
\el
{\bf Proof.} Let $(i,f,j)\in S$. Since $\e$ is perfect, by Result \ref{202}, $\e\circ\delta_r=S\times S$. Hence there is $(m,h,n)\in S$ such that $((i,f,j),( m,h,n))\in\e$ and $((m,h,n), (0,e_G, 0))\in\delta_r$. In view of the group-forgetful homomorphism of $S$ onto $B$, the latter implies that $((m,n), (0,0))\in\delta_r$ in $B$, so $m\leq 0$ whence $m=0$. We have shown that $(0,h,n)\in(i,f,j)\e$ for some $h\in G$ and $n\in\N$.

Now consider $(j, f^{-1}, i)\in S$. By the argument of the preceding paragraph, there exist $h'\in G$ and $m\in\N$ such that $((0, h', m), (j, f^{-1}, i))\in\e$ whence $((0, h', m)^{-1}, (j, f^{-1}, i)^{-1})\in\e$, that is, $((m, (h')^{-1}, 0), (i, f, j))\in\e$. Denoting $(h')^{-1}$ by $g$, we obtain $(m,g,0)\in(i,f,j)\e$. \epr
\bt\label{410} Let $\e=\e_{(N,z,k)}$ with $k\geq 1$. Then $\e$ is perfect if and only if 
\setcounter{equation}{10}
\be (\forall\,n\geq 1)(\forall\,x\in G)\; G\alpha^n\cap Nx\ne\emptyset.
\label{eq:15th}
\ee
\et
{\bf Proof.}  Assume that $\e$ is perfect. Take any $x\in G$ and $n\geq 1$. Then $(0,x,n)\in S$ and, by Lemma \ref{409}, there exists $(m,g,0)\in S$ such that $((m,g,0), (0,x,n))\in\e$. According to (\ref{eq:13th}), there is $l\in\Z$ such that $n+m=kl$ and $x=x\alpha^0\in Nz^l(g\alpha^n)$. Since  $kl=m+n\geq 1$, we conclude that $l\geq 1$. Recall that $Nz^l=N(z^l)\alpha^n$. Therefore $x\in N((z^l)\alpha^n)(g\alpha^n)=N(z^lg)\alpha^n$. It follows that $(z^lg)\alpha^n\in Nx$. Thus $G\alpha^n\cap Nx\ne\emptyset$.   

Conversely, suppose that (\ref{eq:15th}) holds. Let us show that $\e\circ\delta_r=S\times S$. Take arbitrary $(m, g, n), (i, f, j)\in S$. Then $((m, g, n), (i, f, j))\in \e\circ\delta_r$ if there exists $(p, h, q)\in S$ such that $(p, h, q)\in (m, g, n)\e$ and $((p, h, q), (i, f, j))\in\delta_r$. Using the congruence $\zeta_k$ on the bicyclic semigroup $B$ introduced in the first paragraph of this section, we can write (\ref{eq:13th}) as follows:
\be
(p,h,q)\!\in\!(m,g,n)\e\!\Longleftrightarrow\! ((m,n), \!(p,q))\!\!\in\!\zeta_k\text{ and }h\alpha^n\!\!\in\!\!Nz^l(g\alpha^q)\text{ for }l\!=\!\frac{(q\!-\!p)\!-\!(n\!-\!m)}{k},
\label{eq:16th}
\ee 
Recall that $\zeta_k$ is a perfect congruence on $B$. Hence $\zeta_k\circ\delta^B_r=B\times B$ (we denote by $\delta^B_r$ the right divisibility relation on $B$ to distinguish it from the right divisibility relation $\delta_r$ on $S$). It follows that there exists $(p, q)\in B$ such that $((m,n), \!(p,q))\!\in\!\zeta_k$ and $((p, q), (i, j))\in\delta^B_r$. Since (\ref{eq:15th}) holds, it is guaranteed that there is $h\in G$ such that $h\alpha^n\!\in\!Nz^l(g\alpha^q)$ where $l\!=\![(q-p)-(n-m)]/k$. Thus we have found $(p, h, q)\in S$ such that (\ref{eq:16th}) is satisfied.

Since $((p, q), (i, j))\in\delta^B_r$, there is $(s, t)\in B$ such that $(i, j)=(p, q)(s, t)$, and it can be assumed here that $q\leq s$, so that ${\rm min}\{q, s\}=q$; indeed, if $(i, j)=(p, q)(s, t)$ where $q>s$, then we also have $(i, j)=(p, q)(q, q+t-s)$. Let $x=(h\alpha^{s-q})^{-1}f$ whence $f=(h\alpha^{s-q})x$. Since $(i, j)=(p, q)(s, t)$ and $f=(h\alpha^{s-{\rm min}\{q, s\}})(x\alpha^{q-{\rm min}\{q, s\}})$, it follows that $(i, f, j)=(p, h, q)(s, x, t)$ and therefore $((p, h, q), (i, f, j))\in\delta_r$.

We have demonstrated that there exists $(p, h, q)\in S$ such that $((m, g, n), (p, h, q))\in\e$ and $((p, h, q), (i, f, j))\in\delta_r$, so that $((m, g, n), (i, f, j))\in \e\circ\delta_r$. Therefore $\e\circ\delta_r=S\times S$ and hence, by Result \ref{202}, $\e$ is a perfect congruence on $S$. The proof is complete. \epr
\setcounter{theorem}{12}
\bc\label{413}
Let $\e=\e_{(N,z,k)}$ with $k\geq 1$. If $G\alpha=\{e_G\}$ or if $G\alpha=G$, then $\e$ is perfect. 
\ec
{\bf Proof.} Assume that $G\alpha=\{e_G\}$. Since $N$ is $\alpha$-invariant, $G=\{e_G\}\alpha^{-1}\subseteq N\alpha^{-1}=N$, so $N=G$. Hence (\ref{eq:15th}) trivially holds, and $\e$ is perfect by Theorem \ref{410}.\\
\indent Now suppose that $G\alpha=G$. Then $G\alpha^n=G$ for all $n\geq 1$, so (\ref{eq:15th}) is obviously true and again, by Theorem \ref{410}, $\e$ is perfect. \epr
\vspace{0.02in}\\
\indent Corollary \ref{413} shows that if $\alpha$ is a surjective endomorphism of $G$, then $\e_{(N,z,k)}$ (with $k\geq 1$) is perfect but the converse, in general, is not true. However, it might be of interest to note that the converse holds under the additional assumption that $N\alpha=N$.
\bp\label{414}
Let $\e=\e_{(N,z,k)}$ with $k\geq 1$. If $G\alpha=G$, then $\e$ is perfect and $N\alpha=N$. Conversely, if $\e$ is perfect and $N\alpha=N$, then $\alpha$ is a surjective endomorphism of $G$.
\ep
{\bf Proof.} Suppose $G\alpha=G$. Then $\e$ is perfect by Corollary \ref{413}. Let $b\in N$. Then $a\alpha=b$ for some $a\in G$. Now $a\in N\alpha^{-1}=N$ since $N$ is $\alpha$-invariant, so that $b\in N\alpha$. Thus $N\subseteq N\alpha$, and the converse inclusion holds because $N$ is $\alpha$-admissible. Therefore $N=N\alpha$.\\
\indent Conversely, assume that $\e$ is perfect and $N\alpha=N$. Let $h\in G$. Consider $(0, h, 1)\in S$. By Lemma \ref{409}, there is $(m, g, 0)\in S$ such that $((m, g, 0), (0, h, 1))\in\e$. Since $(0, h, 1)\in (m, g, 0)\e$, it follows from (\ref{eq:13th}) that $m+1=kl$ for some $l\geq 1$ and $h=h\alpha^0\in Nz^l(g\alpha)$. Since $N=N\alpha$, we have $h\in  Nz^l(g\alpha)=(N\alpha)(z^l\alpha)(g\alpha)$, so that $h=(a\alpha)(z^l\alpha)(g\alpha)=(az^lg)\alpha$ for some $a\in N$. Thus $\alpha$ is a surjection of $G$ onto $G$. \epr
\vspace{0.01in}\\
\begin{center}
{\bf Acknowledgement}
\end{center}
\sm Preliminary versions of results proved in this paper were obtained in the summer of 2013 by the Algebra Team of the REUT Program at the CSU, Chico. Members of the team were students -- Katherine Grimshaw (North Carolina State University), Anthony Kling (Cal Poly, San Luis Obispo), Freda Li (Bryn Mawr College), and a teacher -- Therese Landry (Chair of the Mathematics Department of the Orthopaedic Hospital Medical Magnet High School, Los Angeles). Their research was supervised and directed by Simon M. Goberstein. The authors are grateful to the National Science Foundation for the financial support and to the administration, faculty and staff of the CSU, Chico for creating excellent working conditions. 
\vspace{0.01in}\\

\end{document}